\documentclass[11pt]{amsart}

\usepackage{color,graphicx,amssymb,latexsym,amsfonts,txfonts,amsmath,amsthm}
\usepackage{pdfsync,bookmark}

\usepackage{hyperref}
\hypersetup{
    colorlinks=true,       % false: boxed links; true: colored links
    linkcolor=blue,          % color of internal links
    citecolor=blue,        % color of links to bibliography
    filecolor=blue,      % color of file links
    urlcolor=blue           % color of external links
}

\begin{document} 
\newcommand{\s}{\vspace{0.2cm}} 

\newtheorem{theo}{Theorem} 
\newtheorem{prop}{Proposition}
\newtheorem{coro}{Corollary}
\newtheorem{lemm}{Lemma}
\newtheorem{claim}{Claim}
\newtheorem{example}{Example}
\theoremstyle{remark}
\newtheorem{rema}{\bf Remark} 
\newtheorem{defi}{\bf Definition}

\title{An effective descent of arithmetical real algebraic varieties}
\date{}

\author{Rub\'en A. Hidalgo} 


\address{Departamento de Matem\'atica y Estad\'{\i}stica, Universidad de La Frontera. Temuco, Chile}
\email{ruben.hidalgo@ufrontera.cl}
\thanks{Partially supported by Project Fondecyt 1150003 and Anillo ACT 1415 PIA-CONICYT}

\begin{abstract}
Let $X$ be a complex smooth algebraic variety admitting a symmetry $L$, that is, an antiholomorphic automorphism of order two. If both, $X$ and $L$ are defined over $\overline{\mathbb Q}$, then Koeck, Lau and Singerman showed the existence of a complex smooth algebraic variety $Z$ admitting a symmetry $T$, both defined over ${\mathbb R} \cap \overline{\mathbb Q}$, and of an isomorphism $R:X \to Z$ so that $R \circ L \circ R^{-1}=T$. The provided proof is existential and, if explicit equations for $X$ and $L$ are given over $\overline{\mathbb Q}$, then it is not described how to get the explicit equations for $Z$ and $T$ over ${\mathbb R} \cap \overline{\mathbb Q}$. In this paper we provide an explicit rational map $R$ defined over ${\mathbb Q}$ so that $Z=R(X)$ is defined over ${\mathbb R} \cap \overline{\mathbb Q}$ and with $T=R \circ L \circ R^{-1}$ being the usual conjugation map.
\end{abstract}

\maketitle

%%%%%%%%%%%%%%%%%%
%%%%%%%%%%%%%%%%%%
\section{Introduction}
Let $X \subset {\mathbb C}^{n}$ be complex smooth algebraic variety and set $\overline{X}=J_{n}(X)$, where $J_{n}(x_{1},\ldots,x_{n})=(\overline{x_{1}},\ldots,\overline{x_{n}})$ is the usual conjugation map on ${\mathbb C}^{n}$. A {\it symmetry} of $X$ is an antiholomorphic automorphism $L:X \to X$ of order two so that $J_{n} \circ L:X \to \overline{X}$ is a biregular isomorphism; in this case, we say that $X$ is {\it symmetric} and that the pair $(X,L)$ is a {\it real algebraic variety}.  Two real algebraic varieties $(X,L)$ and $(Y,T)$ are {\it isomorphic} if there is a birational isomorphism $R:X \to Y$ with $T=R  \circ L \circ R^{-1}$.

A symmetric algebraic variety may have different symmetries, even different conjugacy classes of them (inside its group of automorphisms). In complex dimension one (i.e., Riemann surfaces) the number of (conjugacy classes of) symmetries is well known (see, for instance, \cite{BCGG,BGI}).

Let ${\mathcal Q}$ be a subfield of ${\mathbb C}$ and $(X,L)$ a real algebraic variety. We say that (i) $(X,L)$ is {\it defined} over ${\mathcal Q}$ if both, $X$ and $J_{n} \circ L$, are defined over ${\mathcal Q}$, and (ii) that $(X,L)$ is {\it definable} over ${\mathcal Q}$ if there is a real algebraic variety $(Y,T)$ defined over ${\mathcal Q}$ and isomorphic to $(X,L)$.  As a consequence of Weil's descent theorem \cite{Weil} (see also \cite{Silhol}), every real algebraic variety is definable over ${\mathbb R}$. 

In \cite{K-L,K-S}  Koeck, Lau and Singerman showed that a real algebraic variety $(X,L)$ which is definable over the field $\overline{\mathbb Q}$ of algebraic number is also definable over ${\mathbb R} \cap \overline{\mathbb Q}$. 
More generally, if $(X,L)$ is defined over a subfield ${\mathcal Q}$ of ${\mathbb C}$ so that $J_{1}({\mathcal Q})={\mathcal Q}$ (i.e., ${\mathcal Q}$ is {\it conjugate-invariant}), then it is also definable over ${\mathbb R} \cap {\mathcal Q}$ (Theorem \ref{teo1}).
The proof, based on Weil's descent theorem \cite{Weil}, is existential, that is, given the explicit equations for $X$ and $L$ over ${\mathcal Q}$, there is no explained how to construct explicitly $R$. In this paper, we describe an explicit rational isomorphism map $R$, defined over ${\mathcal Q}$, so that $R:X \to Z=R(X)$ is an isomorphism, $Z$ is defined over ${\mathbb R} \cap {\mathcal Q}$ and with $T=R \circ L \circ R^{-1}$ being the usual conjugation map (see Theorem\ref{teodos}). We hope that this algorithm may be of use in the construction of explicit examples.

%%%%%%%%%%%%%%%%%%%%%%%%%%
%%%%%%%%%%%%%%%%%%%%%%%%%%
\section{Main results}\label{Sec:main}

%%%%%%%%%%%%%%
\subsection{}
Let us start with the following simple application of Weil's descent theorem, a mild generalization of the results in \cite{K-L,K-S}.

\begin{theo}\label{teo1}
Let ${\mathcal Q}$ be a conjugate-invariant subfield of ${\mathbb C}$. If $(X,L)$ is a real algebraic variety defined over ${\mathcal Q}$, then $(X,L)$ is definable over ${\mathcal Q} \cap {\mathbb R}$.
\end{theo}

\begin{rema}
 Assume that ${\mathcal K}$, ${\mathcal N}$ and ${\mathcal L}$ are subfields of ${\mathbb C}$ so that ${\mathcal L}$ contains ${\mathcal K}$ as an algebraically closed subfield (i.e., the only ${\mathcal K}$-algebraic numbers of ${\mathcal L}$ belongs to ${\mathcal K}$) and ${\mathcal N}$ is a finite Galois extension of ${\mathcal K}$. Let $X$ be a complex smooth algebraic variety, which is definable over ${\mathcal L}$ and also over ${\mathcal N}$ (maybe by different models). In \cite{K-L}, as a consequence of Weil's descent theorem, it is shown the existence of an isomorphism $R:X \to Z$, where $Z$ is defined over ${\mathcal K}$ (the given proof does not provide a method to obtain explicitly an isomorphism $R:X \to Z$). We should observe that this result does not implies Theorem \ref{teo1}; for instance  take ${\mathcal Q}={\mathbb Q}(i)$. 
\end{rema}

If the conjugate-invariant subfield ${\mathcal Q}$ is also algebraically closed (for instance, ${\mathcal Q} =\overline{\mathbb Q}$ or ${\mathcal Q}=\overline{{\mathbb Q}(\pi, \overline{\pi})}$), then Theorem \ref{teo1} can be written as follows.

\begin{coro}\label{teomain}
Let ${\mathcal Q}$ be an algebraically closed conjugate-invariant subfield of ${\mathbb C}$.
If $X$ is a symmetric variety definable over ${\mathcal Q}$,
whose group of birational automorphisms is finite, then $X$ is definable over ${\mathcal Q} \cap {\mathbb R}$.
\end{coro}
\begin{proof}
Let $L:X \to X$ be a symmetry of $X$. We claim that $L$ is defined over ${\mathcal Q}$. In fact, if $\eta \in {\rm Gal }({\mathbb C}/{\mathcal Q})$, then we have the symmetry $L^{\eta}:X \to X$. So, there is a birational automorphism $t$ of $X$ so that $L^{\eta}=L \circ t$. 
 Set $K=\{\eta \in {\rm Gal }({\mathbb C}/{\mathcal Q}): L=L^{\eta}\}$ and let $\mathcal U$ be the fixed field of $K$. Assume $L=(L_{1},\ldots,L_{n})$, where $L_{j}$ is a rational map of the form $r_{j}/s_{j}$, where $r_{j}$ and $s_{j}$ are relatively prime polynomials. We may assume the leading coefficient of $s_{j}$ to be equal to $1$.
 The equality $r_{j}^{\eta}/s_{j}^{\eta}=r_{j}/s_{j}$, for $\eta \in K$,  implies that te set of zeroes of $r_{j}$ (and the set of zeroes of $s_{j}$) is invariant under $K$. In particular, $r_{j}^{\eta}=a_{\eta}r_{j}$ and $s_{j}^{\eta}=b_{\eta}s_{j}$, for suitable $a_{\eta}, b_{\eta} \neq 0$. Since $r_{j}^{\eta}/s_{j}^{\eta}=r_{j}/s_{j}$, we also have that $b_{\eta}=a_{\eta}$. 
 As the leading coefficient of $s_{j}$ is equal to $1$, we must have $b_{\eta}=1$; so $r_{j}$ and $s_{j}$ are polynomials defined over ${\mathcal U}$. We conclude that $L$ is defined over ${\mathcal U}$.
Now, as the group of birational automorphisms of 
$X$ is finite, it follows that $K$ is a subgroup of finite index of ${\rm Gal }({\mathbb C}/{\mathcal Q})$; so $\mathcal U$ is a finite extension of ${\mathcal Q}$. As ${\mathcal Q}$ is assumed to be algebraically closed, it follows that ${\mathcal U}={\mathcal Q}$ and we obtain the desired result for $L$. Now, we may apply Theorem \ref{teo1}.
\end{proof}

%%%%%%%%%%%%%%%%%%%%%
\subsection{The explicit construction} 
Let ${\mathcal Q}$ a conjugate-invariant subfield of ${\mathbb C}$ and let 
$(X,L)$ be a real algebraic variety defined over ${\mathcal Q}$. Theorem \ref{teo1} asserts the existence of a real algebraic variety $(Z,T)$, defined over ${\mathcal Q} \cap {\mathbb R}$, and of an isomorphism $R:X \to Z$ so that $T=R \circ L \circ R^{-1}$. Of course, if ${\mathcal Q}$ is a subfield of ${\mathbb R}$, then $R=I$ suffices. So, from now on, we assume that ${\mathcal Q}$ is not a 
subfield of ${\mathbb R}$ and that $(X,L)$ is not already defined over ${\mathcal Q} \cap {\mathbb R}$. Under these assumptions, we proceed to the construction of a rational map $R$ defined over ${\mathcal Q}$ so that: $R:X \to Z=R(X)$ is an isomorphism, $Z$ is defined over ${\mathcal Q} \cap {\mathbb R}$ and $R \circ T \circ R^{-1}$ is the complex conjugation.

Assume $X \subset {\mathbb C}^{n}$ is defined by the polynomials $P_{1}(x_{1},\ldots,x_{n}),\ldots,P_{s}(x_{1},\ldots,x_{n}) \in {\mathcal Q}[x_{1},\ldots,x_{n}]$.

\subsubsection{\bf Step 1: Can assume $X \cap \overline{X}= \emptyset$}
In fact, if $X \cap \overline{X} \neq \emptyset$, then we may change $(X,L)$ by an isomorphic real algebraic variety $(X_{0},L_{0})$ so that $X_{0} \cap \overline{X_{0}} = \emptyset$ as follows. Choose a point $\alpha \in {\mathcal Q}-{\mathbb R}$ and
consider the algebraic variety $X_{0} \subset {\mathbb C}^{n+1}$ (by adding the extra coordinate $x_{n+1}$) defined by the polynomials defining $X$ and the extra polynomial $P_{s+1}(x_{1},\ldots,x_{n+1})=x_{n+1}-\alpha$ (that is, $X_{0}=X \times \{\alpha\}$). The map $H(x_{1},\ldots,x_{n})=(x_{1},\ldots,x_{n},\alpha)$ induces a biregular isomorphism between $X$ and $X_{0}$ (the projection $(x_{1},\ldots,x_{n+1}) \mapsto (x_{1},\ldots,x_{n})$ induces its inverse). It is now clear that $X_{0} \cap \overline{X_{0}}=\emptyset$.
In this new model $X_{0}$ the symmetry is given as 
$$L_{0}=H \circ L \circ H^{-1}:X_{0} \to X_{0}: 
(x_{1},\ldots,x_{n},\alpha) \mapsto (L(x_{1},\ldots,x_{n}),\alpha).$$

\subsubsection{\bf Step 2: Construction of $R$}
By Step 1, we may assume  $X \cap \overline{X} = \emptyset$. 
Consider the following polynomials in ${\mathbb Z}[x_{1},\ldots,x_{n},z_{1},\ldots,z_{n}]$:
\begin{equation}\label{equations1}
\left\{
\begin{array}{lcll}
t_{1,j}&=&x_{j}+z_{j};& j=1,\ldots,n,\\
t_{2,j}&=&x_{j}z_{j};& j=1,\ldots n,\\
t_{k,j}&=&x_{k-2}x_{k-2+j}+z_{k-2}z_{k-2+j};& k=3,\ldots,n, \; j=1,\ldots, n+2-k,\\
t_{n+1,1}&=&x_{n-1}x_{n}+z_{n-1}z_{n}.
\end{array}
\right\}
\end{equation}

For $x=(x_{1},\ldots,x_{n}) \in X$,  set
$z=(z_{1},\ldots,z_{n})=J_{n}(L(x)) \in \overline{X}$ and consider the rational map
$$R:X \to {\mathbb C}^{(n^{2}+3n)/2}:  
x \mapsto t=(t_{1,1},\ldots,t_{n+1,1}).$$

Observe that $R$ is defined over ${\mathcal Q}$. 
Next result states that $R$ as constructed above is the one we are searching.

\begin{theo}\label{teodos}
Under the above assumptions, the following hold.
\begin{enumerate}
\item $Z=R(X)$ is defined over ${\mathcal Q} \cap {\mathbb R}$.

\item $R:X \to Z$ is a birational isomorphism (biregular if $J \circ L$ is a polynomial), defined over ${\mathcal Q}$.

\item $R \circ L \circ R^{-1}=J_{(n^{2}+3n)/2}$.
\end{enumerate}

\end{theo}

%%%%%%%%%%%%%%%%%
%%%%%%%%%%%%%%%%%
\section{Proof of Theorems \ref{teo1} and \ref{teodos}}\label{Sec:Proof}
We only need to care care of the case when the conjugate-invariant subfield ${\mathcal Q}$ is not already a subfield of ${\mathbb R}$.
So ${\mathcal Q}$ is a Galois extension of degree two of ${\mathcal K}={\mathcal Q} \cap {\mathbb R}$.
Let $\Gamma:={\rm Gal}({\mathcal Q}/{\mathcal K})=\langle \sigma(z)=\overline{z}\rangle \cong {\mathbb Z}_{2}$.

Let $P_{1},\ldots,P_{s} \in {\mathcal Q}[x_{1},\ldots,x_{n}]$ be a set of generators of the ideal of $X \subset {\mathbb C}^{n}$ and $L=(L_1,\ldots,L_n)$, where $J_{1} \circ L_j \in {\mathcal Q}(x_1,\ldots,x_n)$.  Observe that $\overline{X}=J_{n}(X)=X^{\sigma}$. If $\sigma_{2}=\sigma$ and $\sigma_{1}=e$ is the identity, then  
we set $f_{2}=f_{\sigma_{2}}=J_{n} \circ L:X \to \overline{X}$, which is an 
isomorphism defined over ${\mathcal Q}$ (since $L$ and $J$ are defined over ${\mathcal Q}$), and  $f_{1}=f_{\sigma_{1}}=I$ the identity.

%%%%%%%%%%%%%%%%%
\subsection{Proof of Theorem \ref{teo1}}
As the collection $\{f_{1}, f_{2}\}$ satisfies Weil's co-cycle conditions $f_{ij}=f_{j}^{\sigma_{i}} \circ f_{i}$, for $i,j \in \{1,2\}$,  it follows from Weil's descent theorem \cite{Weil} the existence of a complex algebraic variety $Z$, defined over ${\mathcal K}$, and an isomorphism $R:X \to Z$, defined over ${\mathcal Q}$, so that $R=R^{\sigma_{j}} \circ f_{\sigma_{j}}$, for $j=1,2$. Now, if $T=R \circ L \circ R^{-1}$, then 
$$T^{\sigma}=R^{\sigma} \circ L^{\sigma} \circ (R^{-1})^{\sigma}=R^{\sigma} \circ L^{\sigma} \circ (R^{\sigma})^{-1}=
R^{\sigma} \circ L^{\sigma} \circ (f_{2} \circ R^{-1})=
R^{\sigma} \circ L^{\sigma} \circ (J_{n} \circ L \circ R^{-1})=$$
$$= R^{\sigma} \circ (J_{n} \circ L \circ J_{n}) \circ J_{n} \circ L \circ R^{-1}= R^{\sigma} \circ (J_{n} \circ L) \circ L \circ R^{-1}=(R^{\sigma} \circ f_2 )\circ L \circ R^{-1}=R \circ L \circ R^{-1}=T,$$
that is, $T$ is defined over ${\mathcal K}$.

%%%%%%%%%%%%%%%%%%%%%%%%%%%%%
\subsection{Proof of Theorem \ref{teodos}}
Now we are assuming that $\overline{X} \cap X=\emptyset$. Let us consider the explicit rational map, defined over ${\mathcal Q}$, 
\begin{equation*}
\Phi:X \subset {\mathbb C}^{n} \to {\mathbb C}^{n} \times {\mathbb C}^{n}={\mathbb C}^{2n}: \quad 
x \mapsto \left(f_{1}(x),f_{2}(x)\right)=\left(x,z \right).
\end{equation*}

As each $f_{j}$ is an isomorphism, $\Phi:X \to \Phi(X)$ is a birational isomorphism whose inverse is given by the restriction to $\Phi(X)$ of the projection map $\pi:{\mathbb C}^{n} \times  {\mathbb C}^{n} \to {\mathbb C}^{n}: (x,z) \mapsto x$.

\begin{rema}
If $f_{2}$ is polynomial (i.e. when $L$ is polynomial), the map $\Phi:X \to \Phi(X)$ is a biregular isomorphism and
$\Phi(X)$ is described by the equations
\begin{equation*}%\label{ecuaciones1-1}
\left\{ \begin{array}{c}
z=f_{2}(x),\\
P_{1}(x)=0, \ldots, P_{s}(x)=0.
\end{array}
\right\}.
\end{equation*}

If $$f_{2}(x)=\left( \frac{r_{1}(x)}{s_{1}(x)},\dots,\frac{r_{n}(x)}{s_{n}(x)}\right),$$
where $r_{j}, s_{j} \in {\mathcal Q}[x_{1},\ldots,x_{n}]$ are relatively prime, then 
$\Phi(X)$ is defined by the equations
\begin{equation*}%\label{ecuaciones1}
\left\{ \begin{array}{c}
z_{1} s_{1}(x)=r_{1}(x), \ldots, z_{n} s_{n}(x)=r_{n}(x),\\
P_{1}(x)=0, \ldots, P_{s}(x)=0,\\
P_{1}^{\sigma_{2}}(z)=0, \ldots, P_{s}^{\sigma_{2}}(z)=0,
\end{array}
\right\}.
\end{equation*}

In any of the above situations, $\Phi(X)$ still defined over ${\mathcal Q}$. 
\end{rema}

\begin{lemm}\label{simetria} 
The symmetry that induces $L$, by $\Phi$, on $\Phi(X)$ is given by
$$L_{\Phi}(x,z)=(\overline{z}, \overline{x}).$$
\end{lemm}
\begin{proof}
$$L_{\Phi}(x,z)=\Phi \circ L \circ \Phi^{-1}(x,z)=\Phi(L (x))=(L(x), J_{n} \circ L (L(x)))=(\overline{z},\overline{x}).$$
\end{proof}

Each $\sigma_{i} \in \Gamma$ produces a permutation $\theta(\sigma_{i}) \in {\mathfrak S}_{2}$ by the rule $\sigma_{i} \sigma_{j}=\sigma_{\theta(\sigma_{i})(j)}$. In fact, $\theta(\sigma_{1})=(1)(2)=e$ is the identity and $\theta(\sigma_{2})=(1,2)$.
We consider the natural isomomorphism (Cayley representation)
\begin{equation*}
\theta:\Gamma \to {\mathfrak S}_{2}: \sigma \mapsto \theta(\sigma).
\end{equation*}

The symmetric group ${\mathfrak S}_{2}$ produces a natural permutation action $\eta({\mathfrak S}_{2})$ on  ${\mathbb C}^{n} \times {\mathbb C}^{n}$ defined as follows. If $x,z \in {\mathbb C}^{n}$, then
\begin{equation*}
\eta(e)(x,z)=(x,z), \quad \eta(1,2)(x,z)=(z,x).
\end{equation*}

In this way, the composition $\widehat{\theta}=\eta \circ \theta$  determines a representation of $\Gamma$ as a subgroup of linear holomorphic automorphisms of ${\mathbb C}^{n} \times {\mathbb C}^{n}$ given by 
\begin{equation*}
\widehat{\theta}(\sigma_{1})(x,z)=\tau_{1}(x,z)=(x,z), \quad
\widehat{\theta}(\sigma_{2})(x,z)=\tau_{2}(x,z)=(z,x)
\end{equation*}
where $x,z \in {\mathbb C}^{n}$.

\begin{lemm} \label{B}
If $j =1,2$, then
\begin{itemize}
\item[(*)] $\widehat{\sigma}_{j}(\Phi(x))=\tau_{j} \circ \Phi \circ f_{\sigma_{j}}^{-1} (\widehat{\sigma}_{j}(x))$.
\item[(**)]  If $y \in \Phi(X)$, then $\widehat{\sigma}_{j}(y) \in \tau_{j}\left(\Phi(X)\right)$. 
\end{itemize}
\end{lemm}
\begin{proof}
Let us recall that, for $\tau, \eta \in \Gamma$ we have the co-cycle relation
$f_{\tau \eta}=f_{\eta}^{\tau} \circ f_{\tau}$.  This co-cycle condition permits to see that, for 
$\eta \in \Gamma$, we have the following sequence of equalities
\begin{equation*} \left\{
\begin{array}{l}
\widehat{\eta}\left(f_{j}(x)\right)=f^{\eta}_{j}\left(\widehat{\eta}(x)\right)=
f^{\eta}_{\sigma_{j}}\left(\widehat{\eta}(x)\right)=
f_{\eta \sigma_{j}}\left(f_{\eta}^{-1}\left(\widehat{\eta}(x)\right)\right)=\\
=f_{\sigma_{\theta(\eta)(j)}}\left(f_{\eta}^{-1}\left(\widehat{\eta}(x)\right)\right)=
f_{\theta(\eta)(j)}\left(f_{\eta}^{-1}\left(\widehat{\eta}(x)\right)\right).
\end{array}\right.
\end{equation*}
\end{proof}

As a consequence of Lemma \ref{B} we have the following  commutative diagram (the top part is just the definition of $\Phi^{\sigma}$ and the bottom part is consequence of the previous lemma)

\begin{equation*}%\label{diagrama0}
\begin{array}[c]{lcl} 
X&\stackrel{\Phi}{\rightarrow}&\Phi(X)\\ 
\downarrow\scriptstyle{\widehat{\sigma}_{j}}&&\downarrow\scriptstyle{\widehat{\sigma}_{j}}\\ 
X^{\sigma_{j}}&\stackrel{\Phi^{\sigma_{j}}}{\rightarrow}&\Phi^{\sigma_{j}}\left(X^{\sigma_{j}}\right)= \Phi(X)^{\sigma_{j}}=\widehat{\sigma}_{j}\left(\Phi(X)\right)=\tau_{j}\left(\Phi(X)\right)\\
\downarrow\scriptstyle{f_{\sigma_{j}}^{-1}}&&\downarrow\scriptstyle{\tau_{j}^{-1}}\\ 
X&\stackrel{\Phi}{\rightarrow}&\Phi(X)
\end{array}
\end{equation*}

\begin{lemm}\label{A}
$\tau_{2}\left(\Phi(X)\right) \cap \Phi(X)=\emptyset$.
\end{lemm}
\begin{proof} Let us assume we have a point $(x,z) \in \tau_{2}\left(\Phi(X)\right) \cap \Phi(X)$. By the definition, $(x,z) \in \Phi(X)$ implies that $x \in X$ and $z=f_{2}(x) \in X^{\sigma_{2}}=\overline{X}$. Since $(x,z) \in \tau_{2}\left(\Phi(X)\right)$, we also have that 
$(z,x) \in \Phi(X)$, that is, $z \in X$. This, in particular, ensures that $X \cap X^{\sigma_{2}}=X \cap \overline{X} \neq \emptyset$, a contradiction.
\end{proof}

Next ingredient in the computational method concerns with the algebra of invariants of a finite group of linear transformations. Let us briefly recall the general facts.
Let ${\mathcal V}$ be a finite dimensional vector space over a field ${\mathcal R}$, say of dimension $n \geq 1$. Let $v_{1},\ldots, v_{n}$ be a basis of ${\mathcal V}^{*}$. The symmetric algebra $S({\mathcal V}^{*})$ over ${\mathcal R}$ can be identified with the algebra ${\mathcal R}[v_{1},\ldots,v_{n}]$. If $G$ is a group acting linearly over ${\mathcal V}$, then it also acts linearly on ${\mathcal V}^{*}$. This provides a linear action on ${\mathcal R}[v_{1},\ldots,v_{n}]$. A theorem due to Hilbert-Noether \cite{Noether1, Noether2} (Chap. 14 in \cite{Procesi}) states that, if $G$ is a finite group, then the algebra of $G$-invariants ${\mathcal R}[{\mathcal V}]^{\Gamma}$ is finitely generated. 

In our situation, $G=\Gamma$ is a cyclic group of order $2$. Noether's bound theorem (see \cite{Noetherbound} for a proof) states that a set of invariant generators can be chosen on the set of polynomials of degree at most $2$. More precisely, the following provides explicit generators.

\begin{lemm}\label{C}
The algebra of $\widehat{\theta}(\Gamma)$-invariant  polynomials is generated by the $(n^{2}+3n)/2$ ones provided in \eqref{equations1} in Section \ref{Sec:main}.
\end{lemm}
\begin{proof}
A polynomial $Q(x_{1},\ldots,x_{n},z_{1},\ldots,z_{n})$ is $\widehat{\theta}(\Gamma)$-invariant if and only if it is symmetric with respect to permutations of the variables $x_{i}$ with $z_{i}$, for all $i=1,\ldots, n$. By Noether's bound theorem a set of invariant generators can be chosen on the set of polynomials of degree at most $2$. 
Now, let us assume we have a polynomial $P$, of degree at most $2$, which is invariant under $\widehat{\theta}(\Gamma)$ and let us consider a monomial of $P$. It will be of the form $x_{i}$ or $z_{j}$ or $x_{i}x_{j}$ or $z_{i}z_{j}$ or $x_{i}z_{j}$. The invariance ensures that its symmetric monomial must also belong to $P$; that is, $z_{i}$ or $x_{j}$ or $z_{i}z_{j}$ or $x_{i}x_{j}$ or $z_{i}x_{j}$, respectively.
As $x_{i}z_{j}+x_{j}z_{i}=(x_{i}+z_{i})(x_{j}+z_{j})-(x_{i}x_{j}+z_{i}z_{j})$ and $x_{i}^{2}+z_{i}^{2}=(x_{i}+z_{i})^{2}-2x_{i}z_{i}$, it follows that $P$ can be generated by polynomials of the form $x_{i}+z_{i}$, $x_{i}z_{i}$ and, for $i \neq j$, $x_{i}x_{j}+z_{i}z_{j}$. 
\end{proof}

Let us consider the polynomial map
\begin{equation*}
\Psi:{\mathbb C}^{n} \times  {\mathbb C}^{n} \to {\mathbb C}^{(n^{2}+3n)/2}; \quad 
\Psi(x,z)=(t_{1,1},\ldots,t_{n+1,1}).
\end{equation*}

\begin{lemm}\label{D}
The map $\Psi$ satisfies the following properties:
\begin{enumerate}
\item $\Psi^{\sigma_{j}}=\Psi$, for every $j=1,2$;
\item $\Psi \circ \tau_{j}=\Psi$, for every $j=1,2$; and
\item $\Psi(w_{1})=\Psi(w_{2})$ if and only if there is some $j \in \{1,2\}$ so that $w_{2}=\tau_{j}(w_{1})$.
\end{enumerate}
\end{lemm}
\begin{proof}
Properties (1) and (2) are trivial to see from the construction of $\Psi$. Now, if $\Psi(x,z)=\Psi(u,v)$, then, for $j=1,\ldots,n$, it holds that
$$x_{j}+z_{j}=u_{j}+v_{j}, \quad x_{j}z_{j}=u_{j}v_{j}.$$

The first equality asserts that $v_{j}=x_{j}+z_{j}-u_{j}$ and now the second one states that
$$u_{j}^{2}-u_{j}(x_{j}+z_{j})+x_{j}z_{j}=0.$$

As the roots of the above are $x_{j}$ and $z_{j}$, we are done.
\end{proof}

Set $V=\Psi({\mathbb C}^{n} \times {\mathbb C}^{n}) \subset {\mathbb C}^{(n^2+3n)/2}$. 
Lemma \ref{D}, together classical invariant theory (see, for instance, \cite{Hilbert,Mumford}), asserts that $V$ is a (singular) algebraic variety of dimension $2n$ whose singular part is $\Psi(\Delta)$, where $\Delta=\{(x,z) \in {\mathbb C}^{2n}: x=z\}$.

\begin{rema}
The algebra of regular maps on $V$ is known to be isomorphic to the algebra ${\mathbb C}[x,z]^{\widehat{\theta}(\Gamma)}$ \cite{CLO}. This can be seen by considering the surjective homomorphism $\xi:{\mathbb C}[t_{1},...,t_{N}] \to {\mathbb C}[x,z]^{\widehat{\theta}(\Gamma)}$ defined by $\xi(p)=p(t_{1},...,t_{N})$, where $N=(n^{2}+3n)/2$.
\end{rema}

The map 
$\Psi:{\mathbb C}^{n} \times  {\mathbb C}^{n} \to V$
is a degree two branched regular holomorphic covering with $\widehat{\theta}(\Gamma)$ as its deck group (see \cite{Tu} for the general concept of a topological branched cover). The branch locus of $\Psi$ is $\Psi(\Delta)$. In particular, $\Psi:{\mathbb C}^{n}\times{\mathbb C}^{n}-\Delta \to V -\Psi(\Delta)$ is a regular holomorphic cover map with $\widehat{\theta}(\Gamma)$ as its deck group.

\begin{lemm}
The map $\Psi:\Phi(X) \to Z=\Psi(\Phi(X))$ defines a biregular isomorphism. In particular, $R=\Psi \circ \Phi:X \to Z$ is defined over 
${\mathcal K}$ and it is a birational isomorphism,  which is biregular if $f_{2}$ is polynomial.
\end{lemm}
\begin{proof}
Set $\Phi(X)=W$. Since $\tau_{2}(W) \cap W)=\emptyset$, we have that the polynomial map $\Psi:W \to Z$ is bijective. Now, set $\widehat{W}=W \cup  \tau_{2}(W)$, which is a reducible affine variety whose irreducible components are $W$ and $\tau_{2}(W)$. Since these irreducible components are pairwise disjoint, we may see that the algebra of regular maps on $\widehat{W}$, say ${\mathcal C}[\widehat{W}]$, is the product of the algebras of regular maps of the components, that is,
$${\mathcal C}[\widehat{W}]= {\mathcal C}[W] \times{\mathcal C}[\tau_{2}(W)].$$

The above isomorphism is just given by the restriction of each regular map of $\widehat{W}$ to each of its irreducible components.
There is also natural isomorphism $\rho:{\mathcal C}[W] \to {\mathcal C}[\widehat{W}]^{\widehat{\theta}(\Gamma)}$, where ${\mathcal C}[\widehat{W}]^{\widehat{\theta}(\Gamma)}$ denotes the sub-algebra of $\widehat{\theta}(\Gamma)$-invariant regular maps on $\widehat{W}$. This isomorphism is given as follows. If $p \in {\mathcal C}[W]$, then $\rho(p)=(p,p \circ \tau_{2}^{-1})$ defines an injective homomorphism. It is clear that every $\widehat{\theta}(\Gamma)$-invariant regular map of $\widehat{W}$ is obtained in that way (so $\rho$ is surjective).  On the other hand, ${\mathcal C}[Z]$ is isomorphic to ${\mathcal C}[\widehat{W}]^{\widehat{\theta}(\Gamma)}$. To see this, one may consider the injective homomorphism $\chi:{\mathcal C}[Z] \to {\mathcal C}[\widehat{W}]^{\widehat{\theta}(\Gamma)}$ defined by $\chi(p)=p \circ \Psi$. Now, to see that $\chi$ is onto, we only need to note that $\rho^{-1}(\chi({\mathcal C}[Z]))$ is a sub-algebra of ${\mathcal C}[W]$, that $W$ is irreducible and that $Z$ has the same dimension as $W$. In this way, $\chi^{-1} \circ \rho$ produces an isomorphism between ${\mathcal C}[W]$ and ${\mathcal C}[Z]$, that is, $\Psi:\Phi(X) \to Z$ is a biregular isomorphism and, in particular, $R=\Psi \circ \Phi:X \to Z$ is a birational isomorphism. As $\Psi:W \to Z$ is biregular isomorphism and $\Phi:X \to W$ is biregular if $f_{2}$ is polynomial, then $R:X \to Z$ turns out to be biregular if $f_{2}$ is polynomial. 
\end{proof}

In order to finish the proof, we only need to show that $Z$ is defined over ${\mathcal K}$ since by Lemma \ref{simetria} and the definition of $\Psi$ we have that $R \circ L \circ R^{-1}$ is just the conjugation map.

\begin{prop}\label{defZ}
$Z$ is defined over ${\mathcal K}$.
\end{prop}
\begin{proof}
If $\eta \in {\rm Aut}({\mathbb C}/{\mathcal Q})$, then, as $X^{\eta}=X$, $\Phi^{\eta}=\Phi$ (since $L$ and $X$ are defined over ${\mathcal Q}$) and $\Psi^{\eta}=\Psi$ (since $\Psi$ is defined over the basic field ${\mathbb Q}$), one has that $R^{\eta}=R$; in particular,  $Z^{\eta}=R(X)^{\eta}=R^{\eta}(X^{\eta})=R(X)=Z$. So, $Z$ is defined over ${\mathcal Q}$. We have that ${\mathcal Q}$ is a degree two extension of ${\mathcal K}$ and ${\rm Gal}({\mathcal Q}/{\mathcal K})$ is generated by $\sigma_{2}$. It follows, from $(**)$ in Lemma \ref{B} and (3) in Lemma \ref{C}, that $\widehat{\sigma_{2}}:Z \to Z$ is a bijection, that is, $Z^{\sigma_{2}}=Z$; so $Z$ is defined over ${\mathcal K}$ as desired.
\end{proof}

%%%%%%%%%%%%%%%%
\subsection{Some final computational remarks}

\subsubsection{Equations for $V$}\label{elimina}

It is not difficult to check, for $1\leq i<j \leq n$ (so $n(n-1)/2$ cases), the following equality
\begin{equation}\label{EV}
t_{1,i}^{2} t_{2,j} - t_{1,i} t_{1,j} t_{i+2,j-i} + t_{1,j}^{2} t_{2,i} - 4 t_{2,i} t_{2,j} + t_{i+2,j-i}^{2}=0.
\end{equation}

It follows that, on $V$, we have that $t_{i+2,j-i}$ is rationally determined by the variables $t_{1,i}$, $t_{2,i}$, $t_{1,j}$ and $t_{2,j}$. In particular, the algebraic variety $V$ is a subvariety of the the algebraic variety $\widehat{V}$ defined by the above $n(n-1)/2$ polynomial equations in \eqref{EV}. In worked examples we have seen that $V=\widehat{V}$. We believe that, in the general case, $V=\widehat{V}$.

%%%%%%%%%%%%%%%%%%
\subsubsection{Finding the inverse of $R$}\label{C1}
In order to find explicitly an inverse of $R:X \to Z$, we need to search for rational expressions for the variables $x_{i}$'s in the variables $t_{j}$'s. In fact, we will see that we only need the variables $t_{1,i}$, $t_{2,i}$ and $t_{3,i}$ (see Section \ref{elimina}).
We first note that  
\begin{equation*}
z_{j}=t_{1,j}-x_{j}, \quad j=1,\ldots,n.
\end{equation*}

If $j=2,\ldots,n$, then 
$$t_{3,j-1}+t_{1,j}x_{1}-t_{1,j}t_{1,1}=x_{1}x_{j}+z_{1}z_{j}+(x_{j}+z_{j})x_{1}-(x_{j}+z_{j})(x_{1}+z_{1})=$$
$$=x_{1}x_{j}-x_{j}z_{1}=(x_{1}-z_{1})x_{j}=(x_{1}-(t_{1,1}-x_{1}))x_{j}=(2x_{1}-t_{1,1})x_{j},$$
that is, 
\begin{equation*}%\label{dos}
x_{j}=\frac{t_{3,j-1}+t_{1,j}(x_{1}-t_{1,1})}{2x_{1}-t_{1,1}}, \quad j=2,\ldots,n.
\end{equation*}

All the above states that each of the variables $x_{2},\ldots, x_{n}$ can be expressed rationally in terms of the variables $t_{1,1},\ldots, t_{3,n-1}$ and $x_{1}$. Since $t_{1,1}x_{1}-t_{2,1}=(x_{1}+z_{1})x_{1}-x_{1}z_{1}=x_{1}^{2}$, we also have
\begin{equation}\label{X12}
x_{1}^{2}=t_{1,1}x_{1}-t_{2,1}.
\end{equation}

Using the polynomial equations of $X$ or the first coordinate polynomial equation of $f_{2}$ (which provides $z_{1}$ in terms of $x_{1},\ldots, x_{n}$) and  equality \eqref{X12} (which permits to pass high powers of $x_{1}$ to a degree $1$ power of it) we may obtain $x_{1}$ in terms of $t_{1,1},\ldots,t_{3,n-1}$. In this way, each of the variables $x_{1},\ldots,x_{n}$ is now rationally expressed in the variables $t_{1,1},\ldots,t_{3,n-1}$ and the inverse map $R^{-1}:Z \to X$ is then obtained.

%%%%%%%%%%%%%%%%%%%%%%
\subsubsection{Elimination of variables for $Z$}
Observe (see Section \ref{elimina}) that, when we are restricted to $X$, the coordinates $t_{2+i,j}$ (where $i<j$) are determined by the coordinates $t_{1,1}$, $t_{1,2}$, ...,  $t_{2,n}$; so $Z$ can be described using only these $2n$ coordinates.

%%%%%%%%%%%%%%%
%%%%%%%%%%%%%%%
\section{An Example}\label{Sec:Example}
Let us consider the real algebraic curve $(X,L)$, defined over ${\mathbb Q}(i)$,  where
\begin{equation*}
X :  \left \{ \begin{array}{lllllll}
\,\,\,\,\,1 & + & x_1^2 & + & x_2^2 & =  & 0\\
 -1  & + & x_1^2 & + & x_3^2 & =  & 0\\
\,\,\, i  & + & x_1^2 & + & x_4^2 & =  & 0\\
\end{array}
\right\}
\subset {\mathbb C}^{4},
\end{equation*}
and
\begin{equation*}
L:{\mathbb C}^{4} \to {\mathbb C}^{4}: 
(x_1,x_2,x_3,x_4) \mapsto (-i \;\overline{x_{1}},-i \; \overline{x_{3}},-i \; \overline{x_{2}},-i \;\overline{x_{4}}).
\end{equation*}

By Theorem \ref{teo1}, $(X,L)$ is definable over ${\mathbb Q}={\mathbb Q}(i) \cap {\mathbb R}$. Theorem \ref{teodos} states that the map
\begin{equation*}\begin{array}{c}
R:X \to Z\\
R(x_{1},x_{2},x_{3},x_{4})=(t_{1}=t_{1,1},\ldots,t_{14}=t_{5,1})\\
\\
t_{1}=(1+i)x_{1},\;
t_{2}=x_{2}+i x_{3},\;
t_{3}=x_{3}+i x_{2}, \;
t_{4}=(1+i)x_{4},\\
t_{5}=i x_{1}^{2}, \;
t_{6}=x_{2}+ix_{3}, \;
t_{7}=i x_{2}x_{3},\;
t_{8}=i x_{4}^{2}, \\
t_{9}=x_{1}x_{2}-x_{1}x_{3}, \;
t_{10}=x_{1}x_{3}-x_{1}x_{2},\;
t_{11}=0,\;
t_{12}=0,\\
t_{13}=x_{2}x_{4}-x_{3}x_{4},\;
t_{14}=x_{3}x_{4}-x_{2}x_{4}.
\end{array}
\end{equation*}
provides an isomorphism between the real algebraic curves $(X,L)$ and $(Z,J_{14})$ with $Z$ defined over ${\mathbb Q}$. It is not difficult to check that 
\begin{equation*}
R^{-1}:Z \to X:(t_{1},\ldots,t_{14}) \mapsto (x_{1},x_{2},x_{3},x_{4})=\left( \frac{t_{1}}{1+i}, \frac{t_{2}-i t_{3}}{2}, \frac{t_{3}-i t_{2}}{2}, \frac{t_{4}}{1+i}\right).
\end{equation*}

The curve $Z$ is defined by the following equations 
\begin{equation*}\begin{array}{c}
t_{14}=-t_{13}=t_{4}(t_{3}-t_{2})/2, \; t_{12}=t_{11}=0, \;
t_{10}=-t_{9}=-t_{1}(t_{2}-t_{3})/2,\\
t_{8}=t_{4}^2/2,\;
t_{7}=(t_{2}^{2}+t_{3}^{2})/4,\;
t_{6}=t_{2},\;
t_{5} =t_{1}^2/2,\\
4+t_{2}^{2}-t_{3}^{2}=0, \; 
t_{1}^{2}+t_{2}t_{3}=0, \;
t_{1}^{2}+t_{4}^{2}-2=0.
\end{array}
\end{equation*}

\s
\noindent
\begin{rema}\label{rutina}
\begin{enumerate}
\item The above also asserts that $X$ is isomorphic to 
\begin{equation*}
Y=\left\{\begin{array}{ccccc}
4+w_{2}^{2}-w_{3}^{2}=0\\
w_{1}^{2}+w_{2}w_{3}=0\\
 w_{1}^{2}+w_{4}^{2}-2=0
\end{array}\right\} \subset {\mathbb C}^{4},
\end{equation*}
the isomorphism given by
\begin{equation*}
\widehat{R}:X \to Y: (x_{1},x_{2},x_{3},x_{4}) \mapsto (t_{1},t_{2},t_{3},t_{4})=(w_{1},w_{2},w_{3},w_{4}).
\end{equation*}

\item The curve $X$ admits the group $H=\langle A_{1},A_{2},A_{3},A_{4}\rangle \cong {\mathbb Z}_{2}^{4}$ as subgroup of conformal automorphisms, where 
\begin{equation*}\begin{array}{c}
A_{1}(x_{1},x_{2},x_{3},x_{4})=(-x_{1},x_{2},x_{3},x_{4}), \;
A_{2}(x_{1},x_{2},x_{3},x_{4})=(-x_{1},-x_{2},x_{3},x_{4}),\\
A_{3}(x_{1},x_{2},x_{3},x_{4})=(x_{1},x_{2},-x_{3},x_{4}), \;
A_{4}(x_{1},x_{2},x_{3},x_{4})=(x_{1},x_{2},x_{3},-x_{4}).
\end{array}
\end{equation*}

In the model $Y$, these correspond to
\begin{equation*}\begin{array}{c}
A_{1}(w_{1},w_{2},w_{3},w_{4})=(-w_{1},w_{2},w_{3},w_{4}), \;
A_{2}(w_{1},w_{2},w_{3},w_{4})=(w_{1},iw_{3},-iw_{2},w_{4}),\\
A_{3}(w_{1},w_{2},w_{3},w_{4})=(w_{1},-iw_{3},iw_{2},w_{4}), \;
A_{4}(w_{1},w_{2},w_{3},w_{4})=(w_{1},w_{2},w_{3},-w_{4}),
\end{array}
\end{equation*}
which are defined over ${\mathbb Q}(i)$. Observe that the minimal field of definition of the pair $(X,H)$ is ${\mathbb Q}(i)$.

\end{enumerate}
\end{rema}

\s

%%%%%%%%%%%%%%%%%%
%%%%%%%%%%%%%%%%%%

\section{MAGMA implementation}\label{magma}
A simple pseudo-routine in MAGMA \cite{MAGMA} to make the computations as in Theorem \ref{teodos}, under the mild assumption that $L$ is polynomial, is the following. Let us assume ${\mathcal Q}={\mathbb Q}(\alpha)<\overline{\mathbb Q}$ and that $q(t) \in {\mathbb Q}[t]$ is the irreducible polynomial associated to $\alpha$. We have each $z_{j}=z_{j}(x_{1},\ldots,x_{n}) \in {\mathcal Q}[x_{1},\ldots,x_{n}]$. 
 
First we state the ambient spaces. 
\begin{enumerate}
\item[$>$] $Q$:=Rationals(\;);
\item[$>$] $P<t>$:=PolynomialRing($Q$);
\item[$>$] $q:=q(t)$;
\item[$>$] $K<r>$:=SplittingField($q$);
\item[$>$]  $A<x_{1},\ldots,x_{n}>$:=AffineSpace($K$,$n$);
\item[$>$] $B<t_{1,1},\ldots,t_{n+1,1}>$:=AffineSpace($K$, $(n^{2}+3n)/2$);
\end{enumerate}

We introduce the algebraic variety $X \subset {\mathbb C}^{n}$
\begin{enumerate}
\item[$>$] $X$:=Scheme($A$,$[P_{1}(x_{1},\ldots,x_{n}),\ldots,P_{s}(x_{1},\ldots,x_{n})]$);
\end{enumerate}

We introduce the rational map $R:X \to Z$
\begin{enumerate}
\item[$>$] $R$:=map$<A->B|[x_{1}+z_{1},\ldots,x_{n-1}*x_{n}+z_{n-1}*z_{n}]>$;
\end{enumerate}

We now ask for equations of the image $Z=R(X)$. 
\begin{enumerate}
\item[$>$] Image($R$);
\item[$>$] $R(X)$;
\end{enumerate}

\begin{rema}
\mbox{}
\begin{enumerate}
\item Note that in the above one should replace the ``$\dots$" by the corresponding data; for instance, if $n=3$, then the line

 ``$A<x_{1},\ldots,x_{n}>$:=AffineSpace($K$,$n$);" 
 
 should be replaced by the line
 
  ``$A<x_{1},x_{2},x_{3}>$:=AffineSpace($K$,$3$);".

\item Usually MAGMA will provide the defining polynomials over ${\mathcal Q} \cap {\mathbb R}$, but in case that some of these polynomials, say $M \in {\mathcal Q}[t_{1,1},\ldots,t_{n+1,1}]$, is not defined over the desired field, then
we may replace it by the two new trace polynomials ${\rm Tr}(M)=M+\overline{M}, {\rm Tr}(a M)=a M+\overline{a}\overline{M} \in {\mathcal Q} \cap {\mathbb R}[t_{1,1},\ldots,t_{n+1,1}]$, where $\{1,a\}$ is a basis of ${\mathcal Q}$ as a ${\mathcal Q} \cap {\mathbb R}$-vector space and $\overline{M}$ is obtained from $M$ after conjugating its coefficients. Since $M=\lambda_{1} {\rm Tr}(M)+\lambda_{2} {\rm Tr}(a M)$, where $\lambda_{1}= \overline{a}/(\overline{a}-a)$ and $\lambda_{2}=1/(a-\overline{a})$,
 the new set of polynomials, all of them defined over ${\mathcal Q} \cap {\mathbb R}$, will generate the ideal of $Z$. 

\item In the example provided in Section \ref{Sec:Example}, the MAGMA pseudo-routine described above with $q=t^{2}+1$ provides the following 
equations for $Z$ (all of them already over ${\mathbb Q}$):
\begin{equation*}
\left\{\begin{array}{ccccc}
t_{13} + t_{14}=0, \; 
t_{12}=0, \;
t_{11}=0,\;
t_{9} + t_{10}=0,\\
t_{1} + \frac{1}{4} t_{3}t_{10}^{3} - \frac{3}{2}t_{3}t_{10} + \frac{1}{8}t_{4}t_{10}^{3}t_{14} + \frac{1}{8}t_{4}t_{10}t_{14}^{3} +
\frac{3}{4}t_{4}t_{10}t_{14}=0,\\
t_{7} - t_{8} - \frac{1}{2}t_{10}^{2} - \frac{1}{2}t_{14}^{2} + 1=0,\\
t_{8}t_{14}^{2} - t_{8} + \frac{1}{4}t_{10}^{2}t_{14}^{2} + \frac{1}{4}t_{14}^{4} - t_{14}^{2}=0,\\
t_{8}t_{10}^{2} + t_{8} - \frac{1}{4}t_{10}^{2}t_{14}^{2} - \frac{1}{4}t_{14}^{4}=0,\\
t_{10}^{4} + 2t_{10}^{2}t_{14}^{2} - 4t_{10}^{2} + t_{14}^{4} - 4=0,\\
t_{6} - t_{8} - \frac{1}{2}t_{10}^{2} - \frac{1}{2}t_{14}^{2} + 1=0,\\
t_{5} + t_{8} - 1=0,\; 
t_{4}^{2} - 2t_{8}=0,\\
t_{3}t_{14} - \frac{1}{2}t_{4}t_{10}^{2} - \frac{1}{2}t_{4}t_{14}^{2} - t_{4}=0,\\
t_{3}t_{8} - t_{4}t_{8}t_{14} - \frac{1}{4}t_{4}t_{10}^{2}t_{14} - \frac{1}{4}t_{4}t_{14}^{3} + \frac{1}{2}t_{4}t_{14}=0,\\
t_{2} + \frac{1}{2}t_{3}t_{10}^{2} - 2t_{3} + \frac{1}{4}t_{4}t_{10}^{2}t_{14} + \frac{1}{4}t_{4}t_{14}^{3} + \frac{3}{2}t_{4}t_{14}=0,\\
t_{3}t_{4} - 2t_{8}t_{14} - \frac{1}{2}t_{10}^{2}t_{14} - \frac{1}{2}t_{14}^{3} + t_{14}=0,\\
t_{3}^{2} - 2t_{8} - t_{10}^{2} - t_{14}^{2}=0.
\end{array}\right\}.
\end{equation*}
\end{enumerate}
\end{rema}

%%%%%%%%%%%%%%%%

\end{document}